\documentclass[12pt]{amsart}
\usepackage{amssymb}
\usepackage{graphicx}
\usepackage{amsmath}
\usepackage{amsfonts}
\usepackage{subfigure}
\usepackage{hyperref}
\newcommand{\rsphere}{ {\mathbb C_\infty} }

\newcommand{\bd}{\partial}

\newcommand{\sm}{\setminus}

\newtheorem{thm}{Theorem}%[section]

\newtheorem*{thma}{Theorem}
\theoremstyle{definition}
\newtheorem{ques}{Question}
\newtheorem{defn}[thm]{Definition}
\newtheorem*{mconj}{Makienko's Conjecture}

\theoremstyle{remark}

\begin{document}

\title[Buried Points in Julia Sets]{Buried Points in Julia Sets}

\author[C.~P.~Curry]{Clinton P.~Curry}
\email[Clinton P.~Curry]{clintonc@uab.edu}

\author[J.~C.~Mayer]{John C.~Mayer}
\email[John C.~Mayer]{mayer@math.uab.edu}

\address[Clinton P.~Curry and John C.~Mayer]
{Department of Mathematics\\ University of Alabama at Birmingham\\
 Birmingham, AL 35294-1170}



\thanks{We thank Murat Tuncali for inspiration.}

\begin{abstract}
  An introduction to buried points in Julia sets and a list of
  questions about buried points, written to encourage aficionados of
  topology and dynamics to work on these questions.
\end{abstract}

\date{\today}

\maketitle

\begin{center} Dedicated to Bob Devaney on the occasion of his 60th
birthday. \end{center}

\begin{section}{Introduction}
  Let $R:\rsphere\to\rsphere$ be a rational function, where $\rsphere$
  denotes the Riemann sphere.  The \emph{Fatou set} of $R$, denoted
  $F(R)$, is the domain of normality for the family of functions
  $\{R^i \, \mid \, i \in \mathbb N\}$.  A component of the Fatou set
  is called a \emph{Fatou component}.  The \emph{Julia set} of $R$,
  denoted $J(R)$, is the complement of $F(R)$.  General references for
  Juia sets of rational functions are \cite{Beardon:1991},
  \cite{Milnor:2006fr}, and \cite{Carleson:1993vl}; we present a few
  facts below.

  In the case that the degree of $R$ is at least two, the Julia set is
  a non-empty, compact, perfect subset of $\rsphere$.  It is either
  nowhere dense in $\rsphere$ or equal to $\rsphere$.  (We are more
  interested in the former case, and assume it is so henceforth.)  It
  is well-known that $J(R^n)=J(R)$ and $F(R^n)=F(R)$ for any integer
  $n \ge 1$.  The Julia set and Fatou set are each \emph{fully
    invariant} under $R$, meaning that $R^{-1}(J(R)) = J(R)$ and
  $R^{-1}(F(R))=F(R)$.  The restriction $R|_{J(R)}$ is
  \emph{topologically exact}: if $U \subset J(R)$ is open in $J(R)$,
  there exists $n \in \mathbb N$ such that $R^n(U) = J(R)$.

  The notion of a {\em buried point} is purely topological, but we
  consider it here just in the case of Julia sets in $\rsphere$.

  \begin{defn}[Buried Points]
    A point of a Julia set $J(R)$ is said to be \emph{buried} if it does
    not belong to the boundary of a Fatou component.  The set of all
    buried points of $J(R)$  is called the
    \emph{residual Julia set}, denoted $J'(R)$.
  \end{defn}

  Polynomial Julia sets have no buried points, because the Fatou
  component containing $\infty$ has the Julia set as its boundary.
  Curt McMullen \cite{McMullen:1988} presented the first examples of
  rational maps with non-empty residual Julia sets.  He showed that
  functions of the form $z \mapsto z^n+\lambda/z^d$ have Julia sets
  homeomorphic to the product of a Cantor set and a circle whenever
  $1/n+1/d<1$ and $\lambda$ is sufficiently small.  In this case, the
  Julia set is not connected and there are uncountably many components
  of the Julia set which do not intersect the boundary of any Fatou
  component.  Later, John Milnor and Tan Lei \cite{Milnor:1993} and
  Bob Devaney, Dan Look, and David Uminsky \cite{Devaney:2005} exhibited rational
  functions with Julia sets homeomorphic to the Sierpinski carpet (see
  Figure~\ref{fig:carpet}).

  \begin{figure}
    \subfigure[The Sierpinski
    carpet.]{ \label{fig:carpet} \includegraphics[height=2in]{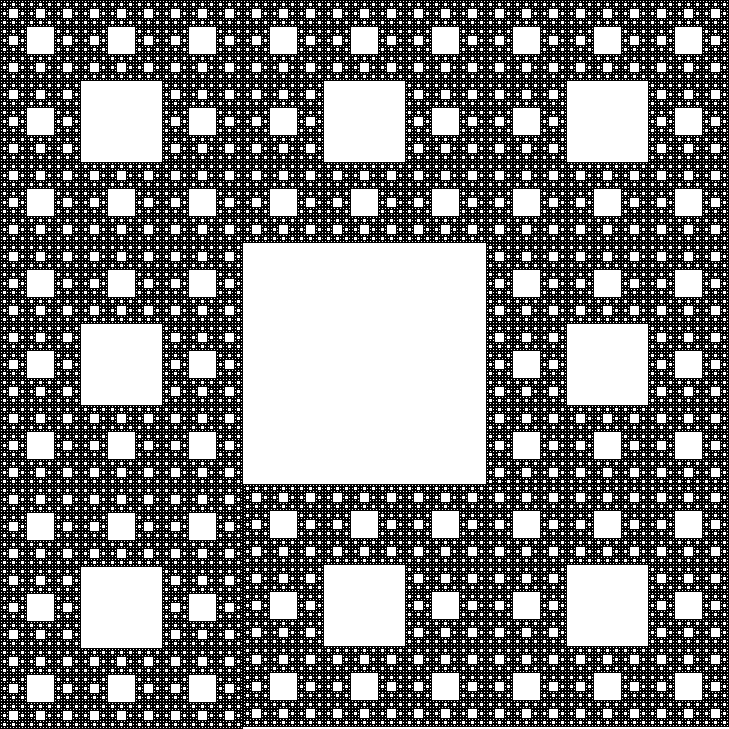}}
    \subfigure[The Sierpsinski
    gasket.]{ \label{fig:triangle} \includegraphics[height=2in]{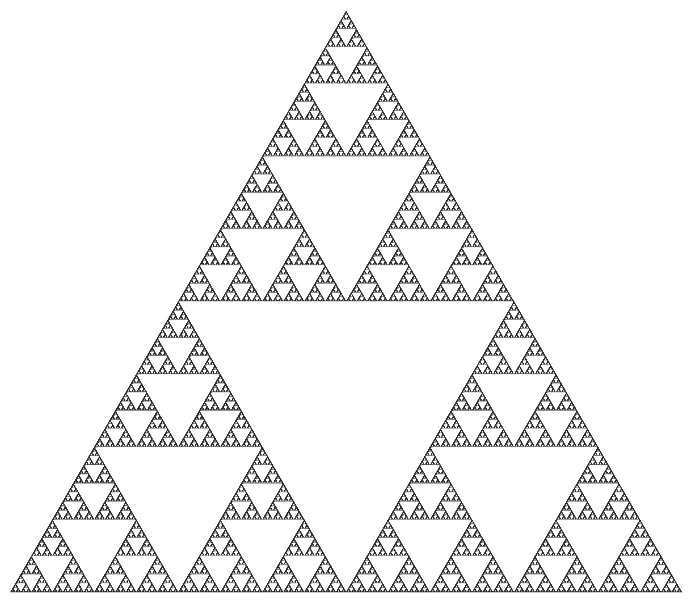}}
  \end{figure}

  Since $R|_{J(R)}$ is topologically exact, we have that the residual
  Julia set is non-empty if and only if the boundary of each Fatou
  component is nowhere dense in $J(R)$.  As a consequence, the
  residual Julia set is a dense $G_\delta$ subset of $J(R)$ whenever
  it is non-empty, as there are only countably many Fatou components.
  It is also nowhere locally compact (because the union of the
  boundaries of Fatou components is also dense).  The residual Julia
  set and the union of boundaries of Fatou components are each fully
  invariant subsets of $J(R)$.

  Interestingly, the known examples of Julia sets which have non-empty
  residual Julia sets are precisely the examples for which no Fatou
  component has a finite grand orbit.  Peter M.~Makienko made the
  following conjecture.

  \begin{mconj}
    Let $R:\rsphere \rightarrow \rsphere$ be a rational function.  The
    Julia set $J(R)$ has buried points if and only if there is no
    completely invariant component of the Fatou set of $R^2$.
  \end{mconj}

  This conjecture was formulated as a possible entry in Sullivan's
  dictionary of correspondences between Kleinian group actions and iteration theory of rational functions in 1990 \cite{Eremenko:1990} as a parallel to a theorem
  of Abikoff.  Note that one direction is easily proved with facts
  already given.  Specifically, if a Fatou component $F$ is completely
  invariant under $R^2$, then $\partial F \subset J(R)$ is also
  completely invariant under $R^2$ and closed, so $\partial F = J(R^2)
  = J(R)$ \cite[Corollary 4.13]{Milnor:2006fr}.  The example $z
  \mapsto \frac 1 {z^2}$ illustrates why one must examine the Fatou
  set of $R^2$.  However, it is known \cite[Theorem
  9.4.3]{Beardon:1991} that a rational map may have at most two
  completely invariant Fatou components, in which case the Julia set
  is a simple closed curve.

  Makienko's conjecture has received attention in the past, with
  results being limited by topological considerations: see Morosawa
  \cite{Morosawa:1997,Morosawa:2000}, Qiao \cite{Qiao:1997lr}, and Sun
  and Yang \cite{Sun:2003}.  Clinton Curry, John Mayer, Jonathan
  Meddaugh, and Jim Rogers \cite{Curry:2008} recently proved the
  following.

  \begin{thma}[Makienko's Conjecture Holds for Decomposable Julia Sets of
    Rational Maps]\label{indecMak}
    If $R$ is a rational function such that $J(R)$ has no buried
    points and $F(R^2)$ has no completely invariant components, then
    $J(R)$ is an indecomposable continuum.
  \end{thma}

  Recall that a continuum is \emph{decomposable} if it can be written
  as the union of two proper subcontinua; otherwise it is
  \emph{indecomposable}.  There are no known examples of Julia sets
  which are indecomposable continua.  In fact, whether or not there
  exists a rational function with an indecomposable continuum as its
  Julia set is a well-known unsolved problem \cite{Mayer:1993}.
\end{section}

\begin{section}{Questions}\label{questions}

  While Makienko's Conjecture is probably the primary motivation for
  dynamicists to consider questions about buried points in Julia sets,
  there are a host of topological questions that one could consider
  that seem interesting from a topological viewpoint.  We present some
  of them here.

  We remarked earlier that the Sierpinski carpet (the Sierpinski
  universal plane curve \cite[1.11]{Nadler:1992ek}, see Figure~\ref{fig:carpet}) appears as a Julia
  set for some rational functions. Sierpinski carpet Julia sets are
  ubiquitous in the families $z \to z^n + \frac{\lambda}{z^d}$ for
  $n\ge 1,\ d\ge 2$ studied by Devaney, Paul Blanchard,  and their students
  \cite{Blanchard:2006rz}.  The set of buried points of the Sierpinski
  carpet is the set of so-called {\em irrational points} of the
  Sierpinski carpet, studied by Krasinkiewicz
  \cite{Krasinkiewicz:1969rm}.  They are one of the two orbits in the
  Sierpinski carpet under its group of homeomorphisms, the points of
  boundaries of complementary components forming the other orbit.
  This set of buried points is one-di\-men\-sion\-al, connected, and locally
  path connected, and it contains a copy of every planar curve
  (one-di\-men\-sion\-al plane continuum). The latter is proved using Whyburn's
  characterization of the Sierpinski carpet \cite{whyburn:1958}.

  The Sierpinski gasket (triangular Sierpinski curve, see
  Figure~\ref{fig:triangle}) and generalizations of it also appear as
  Julia sets in the aforementioned families as studied by Devaney and
  Monica Morena-Rocha \cite{Devaney:2007zl}.  In this case, the set of
  buried points is $0$-di\-men\-sion\-al, and in fact can be shown to be
  homeomorphic to the irrational points on the real line.

  As previously noted \cite{McMullen:1988}, McMullen showed that the set of buried points
  of a Julia set could be the ``irrational'' factors of a Cantor set
  cross a circle.  Such Julia sets occur in all the families $z\to z^n + \frac1{z^n}$
  for $n\ge 3$ and small $\lambda$ \cite{Devaney:2005mc}.

  In \cite{Blanchard:2008,Blanchard:2007}, Blanchard, Devaney, Antonio
  Gajiro, Sebastian Marotta, and Elizabeth Russell study the family $z\to
  z^n+c+\frac{\lambda}{z^n}$ and observe some more complicated residual
  Julia sets.  However,  the residual Julia sets in their examples appear
  to decompose into mutually dense subsets  homeomorphic to one of the
  ``irrational" spaces previously noted.

  \begin{defn}
  A space $X$ is {\em homogeneous} iff there is exactly one orbit in its
  group of homeomorphisms.  That is, given $x,y\in X$, there is a
  homeomorphism $h:X\to X$ with $h(x)=y$.  A space $X$ is {\em
  $1/n$-homogeneous} iff there are exactly $n$ orbits in its group of
  homeomorphisms, each of which is dense.\footnote{The requirement that each
  orbit be dense in $X$ is not part of the standard definition.}
  \end{defn}

 Note that homogeneous is $1/1$-homogeneous. The Sierpinski carpet is an example of a $1/2$-homogeneous space. In each of the aforementioned examples of residual Julia sets, the set of buried points is either a  homogeneous  space or a $1/2$-homogeneous space. That is not to say that the homeomorphisms can necessarily be extended to the Julia set.

  \begin{ques}
    Is there a residual Julia set whose components are not either
    (irrational) points, (irrational) circles, or the irrational
    points of the Sierpinski curve.
  \end{ques}

  \begin{ques}\label{homogeneous}
    Is the residual Julia set always $1/n$-homogeneous for some $n$?  Can the
    homeomorphism be extended to the Julia set?
  \end{ques}

  The reader can readily verify that the answer to the second part of Question~\ref{homogeneous} is
  ``no'' in the case of the Sierpinski gaskets, and ``yes'' when the
  Julia set is itself the Sierpinski carpet.  One can define both stronger and weaker homogeneity properties, which have received attention in topology, and formulate questions similar to Question~\ref{homogeneous}.

  As noted earlier, Krasinkiewicz showed that in the Sierpinski carpet
  any buried point can be carried by a homeomorphism of the carpet to
  any other buried point \cite{Krasinkiewicz:1969rm}.  In
  \cite{Childers:2006fk} it is proved that any buried arc (homeomorphic image of $[0,1]$) can be
  mapped to any other buried arc by a homeomorphism of the carpet.
  This motivates the following question.

  \begin{ques}
    Let $S \subset \rsphere$ be a Sierpinski carpet.  Suppose that
    $X$ and $Y$ are continua in the buried points of $S$ that are equivalently
    embedded in $\rsphere$ (i.e., there
    exists a homeomorphism $h:\rsphere \rightarrow \rsphere$ such that
    $h(X) = Y$).  Is there a homeomorphism $g:S \rightarrow S$ such
    that $g(X) = Y$?
  \end{ques}

  The assumption that $X$ and $Y$ be equivalently embedded in $\rsphere$ is
  necessary, since any homeomorphism from $S$ to itself can be
  extended to a homeomorphism of $\rsphere$ \cite{Krasinkiewicz:1969rm}.

  \begin{ques}
    What does local connectivity imply about the residual Julia set?
    For example, is a Julia set with connected buried point set
    locally connected iff the set of buried points is locally
    connected?
  \end{ques}

  For Julia sets, local connectedness implies connectedness.  This
  follows from the fact that $R|_{J(R)}$ is topologically exact; if $U
  \subset J(R)$ is  connected and open, then the iterate $R^n(U)$
  which equals $J$ must also be connected.  However, there are
  connected polynomial Julia sets which are not locally connected
  \cite[Corollary  18.6]{Milnor:2006fr}.  Recently, Pascale Roesch showed that
  there are ``genuine"  rational Julia sets that are connected but not
  locally connected \cite{Roesch:2006}.  The rational map is not conjugate
  to a polynomial, but the non-local-connectivity is achieved through
  reference to polynomials.  It would be interesting to determine what the
  residual Julia set is in Roesch's examples.

  Suppose that $J(R)$ is a locally connected Julia set with totally disconnected residual Julia set $J'(R)$.
  We conjecture that $J'(R)$ is $0$-di\-men\-sion\-al,
  which would imply that it is homeomorphic to the irrational real
  numbers. However, we have no proof.
  In the case that $J(R)$ is connected, but not locally connected, the situation is
  even less clear.

  \begin{defn}[$0$-Dimensional, Almost $0$-Dimensional]
    A space $X$ is \emph{$0$-di\-men\-sion\-al} at a point $x\in X$ if there is a basis at
     $x$ of open sets
    whose boundaries are empty.  $X$ is \emph{$1$-di\-men\-sion\-al} at $x$ if there is
    a basis at $x$ of open sets whose boundaries are
    $0$-di\-men\-sion\-al.  $X$ is {\em almost $0$-di\-men\-sion\-al} at a point
    $x\in X$ iff $X$ is $1$-di\-men\-sion\-al at $x$ and $x$ has a basis of neighborhoods
    $U$ such that $X\sm U$ is the union of countably many sets with empty
    boundary.
  \end{defn}

  \begin{ques}
    What can be said about the topological dimension of the set of
    buried points?  For example, if the set of buried points is
    totally disconnected, is it $0$-di\-men\-sion\-al?  If not
    $0$-di\-men\-sion\-al, is it almost $0$-di\-men\-sion\-al?
  \end{ques}

  Totally disconnected compact Hausdorff spaces are always
  $0$-di\-men\-sion\-al.  However, it is well-known that totally
  disconnected complete metric spaces can be any dimension.  In
  particular, topologically complete totally disconnected subsets of
  $\rsphere$ can be $1$-di\-men\-sion\-al and not almost $0$-di\-men\-sion\-al.
  Of course, such examples are nowhere locally compact.

  \begin{ques}
    What can be said about the Hausdorff dimension of the set of
    buried points?  For example, can a Julia set and its (nonempty)
    residual Julia set have different Hausdorff dimensions?  In
    particular, what about a Sierpinski gasket Julia set?
  \end{ques}

  Useful references for Hausdorff dimension and related topics include
  \cite{Falconer:2003kx} and \cite{Falconer:1986yq}.

  \begin{ques}
    What does the nature of the post-critical set say about the buried
    point set?
  \end{ques}

  It is known that if the post-critical set is finite, then the Julia
  set is locally connected \cite[Theorem 19.7]{Milnor:2006fr}.

  \begin{ques}
    For a connected Julia set, if a critical value lies on the
    boundary of a Fatou component, is the set of buried points
    disconnected?  When is the set of buried points totally
    disconnected?
  \end{ques}

  Sierpinski gasket Julia sets in the family $z\to
  z^n+\frac{\lambda}{z^d}$ arise when a critical point is on the
  boundary of the Fatou component $B_\infty$ containing $\infty$ and
  the critical value is pre-periodic \cite[Theorem
  3.1]{Devaney:2007zl}.  Since all Fatou components are pre-images of
  $B_\infty$, this serves to (totally) disconnect the residual Julia
  set.

  \begin{ques}
    What can be said about invariant or periodic continua in the set
    of buried points?
  \end{ques}

  For example, in the case of the Sierpinski gasket Julia sets, since boundaries of Fatou components eventually map to $\partial B_\infty$, which is invariant, nearly all of the periodic points {\bf are} buried.

  \begin{ques}
    What can be said about wandering non-degenerate continua in the
    set of buried points?
  \end{ques}

  Blokh and Levin \cite{Blokh:2002yg} bound the number and valence of
  wandering continua for polynomials.

  \begin{ques}
    Assuming the existence of an indecomposable rational Julia set
    with buried points, what can be said about its residual Julia set?
  \end{ques}

  Given a point $x$ in the continuum $X$, the {\em composant} of $x$
  is the union of all proper subcontinua of $X$ containing
  $x$. Composants are dense \cite[5.20]{Nadler:1992ek} in $X$,
  and for an indecomposable continuum are pairwise disjoint
  \cite[Theorem 11.17]{Nadler:1992ek} and uncountable in number
  \cite[Theorem 11.15]{Nadler:1992ek}. There is a relationship between
  buried points and {\em internal composants} of an indecomposable
  plane continuuum studied by Krasinkiewicz
  \cite{Krasinkiewicz:1972lq}.  For example, the boundary of a Fatou
  component is always in an {\em accessible}, hence {\em external}
  composant.  If the residual Julia set is non-empty, then all
  internal composants are in the residual Julia set.  Because
  composants are dense, it seems that part of each external composant
  must also be in the residual Julia set, unless the union of the
  boundaries of infinitely many Fatou components is connected.

\end{section}

\begin{section}{Extensions}

  In view of Theorem~\ref{indecMak}, any counterexample to Makienko's
  conjecture must be exceedingly complicated.  In fact, it is not
  known if the Julia set of a rational function can be as complicated
  as required.  We restate here the questions which appear in
  \cite{Curry:2008,Childers:2006fk, Childers:2006lr, Sun:2003,
    Rogers:1998, Mayer:1993}.
  \begin{ques}
    Can the Julia set of a rational function be an indecomposable
    continuum?
  \end{ques}

\begin{ques}
    Can the Julia set of a rational function contain an
    indecomposable subcontinuum with interior?
  \end{ques}

  Theorems in \cite{Childers:2006fk, Childers:2006lr} indicate that
  the answers to these two questions are the same for polynomial Julia
  sets and for rational functions whose Julia set contains no buried
  points.

  Several authors have extended the study of Makienko's conjecture to
  transcendental entire and meromorphic functions.  Dom\'inguez and
  Fagella \cite{Dominguez:ta} survey the current state of affairs in
  this direction.  Also see Ng, Zheng, and Choi \cite{Ng:2006lr}, who
  prove Makienko's conjecture for locally connected Julia sets of
  certain meromorphic functions.  The techniques of proof of
  Theorem~\ref{indecMak} fail dramatically for trancendental
  functions.  A particular point of difficulty is summarized in the
  following question.
  \begin{ques}
    Let $f$ be a transcendental entire or meromorphic function.  If
    $V$ is a Fatou component such that $\bd V$ is nowhere dense in
    $J(f)$, does it follow that $\bd f(V)$ is also nowhere dense in
    $J(f)$?
  \end{ques}

  It is also worth noting that the proof of Theorem~\ref{indecMak}
  makes use of the non-existence of wandering Fatou components for
  rational functions proved by Dennis Sullivan \cite[Theorem 8.1.2]{Beardon:1991},
  which in general is not true for transcendental functions \cite{Bergweiler:1993}.

  Variations of most of the questions in Section~\ref{questions} could be
  formulated for trancendental functions.  Julia sets of the family
  $z\to \lambda e^z$ would be a place to start.

\end{section}

\bibliographystyle{annotate}

\end{document}